\documentclass{elsart}

\usepackage{amsfonts}
\usepackage{amsmath}
\usepackage[xcolor=pst]{pstricks}

\begin{document}

\begin{frontmatter}
\title{Set maps, umbral calculus, and the chromatic polynomial}

\author{Gus Wiseman\thanksref{vigre}}

\ead{gus@math.ucdavis.edu}

\address{Department of Mathematics,
University of California,
One Shields Ave.,
Davis, CA 95616}

\thanks[vigre]{Research supported in part by NSF VIGRE Grant No. DMS-0135345.}

\begin{abstract}
Some important properties of the chromatic polynomial also hold for any polynomial set map satisfying
\[
p_S(x+y)=\sum_{T\uplus U=S}p_T(x)p_U(y).
\]
Using umbral calculus, we give a formula for the expansion of such a set map in terms of any polynomial sequence of binomial type. This leads to some new expansions of the chromatic polynomial. We also describe a set map generalization of Abel polynomials.
\end{abstract}
\begin{keyword}
chromatic polynomial \sep umbral calculus \sep Abel polynomials
\renewcommand{\MSC}{{\par\leavevmode\hbox{\it2000 MSC:\ }}} 
\MSC 05C15 \sep 05A40 \sep 05A18
\end{keyword}

\end{frontmatter}

\bibliographystyle{hamsplain}

\section{Introduction}\label{section:intro}

In order to motivate our results we will begin by describing two classical expansions of the chromatic polynomial. Let $G$ be a simple graph with finite vertex set $V$. Let $\chi_G(x)$ be the number of proper colorings of $G$ with $x$ colors (assignments of colors $1, 2, \ldots, x$ to the vertices of $G$ so that no two adjacent vertices have the same color). If $e$ is an edge of $G$, let $G\backslash e$ be $G$ with $e$ removed, and let $G/e$ be $G$ with $e$ contracted to a single vertex. Then $\chi_G(x)$ satisfies the \emph{deletion-contraction property},
\[
\chi_G(x)=\chi_{G\backslash e}(x)-\chi_{G/e}(x).
\]
This property can be used to obtain the \emph{subgraph expansion},
\begin{equation}\label{eq:subgraphexpansion}
\chi_G(x)=\sum_{E(H)\subseteq E(G)}(-1)^{|E(H)|}x^{c(H)},
\end{equation}
where the sum is over all spanning subgraphs of $G$ (which correspond to subsets of the edge set $E(G)$), and $c(H)$ is the number of connected components of $H$ (see~\cite[Thm. 10.4]{MR95h:05105} for a proof). In particular, $\chi_G(x)$ is a polynomial, called the \emph{chromatic polynomial} of $G$. The coefficient of $x$ is
\[
\chi^\prime_G(0)=\sum_{\substack{E(H)\subseteq E(G)\\\textrm{connected}}}(-1)^{|E(H)|},
\]
where the sum is over all connected spanning subgraphs. Choosing a spanning subgraph of $G$ is equivalent to choosing a set partition of the vertices and a connected spanning subgraph of the restriction of $G$ to each block. Hence~\eqref{eq:subgraphexpansion} may be expressed in the form
\begin{equation}\label{eq:subgraphexpansion2}
\chi_G(x)=\sum_{\sigma\vdash V}x^{\ell(\sigma)}
\prod_{T\in\sigma}\chi^\prime_{G|_T}(0),
\end{equation}
where the sum is over all set partitions of the vertex set $V$, $\ell(\sigma)$ denotes the number of blocks of $\sigma$, and $G|_T$ is the subgraph induced by $T\subseteq V$.

The second expansion is simpler. A proper coloring of $G$ induces a set partition $\sigma$ of the vertices where two vertices are in the same block if and only if they have the same color. This set partition will be stable, meaning $G$ has no edges with both ends in the same block. Corresponding to such a set partition are $(x)_{\ell(\sigma)}:=x(x-1)(x-2)\cdots(x-\ell(\sigma)+1)$ proper colorings. Hence
\[
\chi_G(x)=\sum_{\substack{\sigma\vdash V\\\textrm{stable}}}(x)_{\ell(\sigma)}.
\]
Since
\begin{displaymath}
\chi_G(1)=\left\{
\begin{array}{ll}1 & \textrm{if $G$ has no edges}\\0 & \textrm{otherwise}\end{array} \right.,
\end{displaymath}
we have
\begin{equation}\label{eq:stablesetexpansion2}
\chi_G(x)=\sum_{\sigma\vdash V}(x)_{\ell(\sigma)}
\prod_{T\in\sigma}\chi_{G|_T}(1).
\end{equation}

Notice the structural similarity between~\eqref{eq:subgraphexpansion2} and~\eqref{eq:stablesetexpansion2}. The present work was motivated by the observation that both expansions are implied by the fact that the chromatic polynomial satisfies the identity
\begin{equation}\label{eq:btype}
\chi_{G}(x+y)=\sum_{T\uplus U=V}\chi_{G|_T}(x)\chi_{G|_U}(y),
\end{equation}
where the sum is over all ordered pairs $(T,U)$ of disjoint sets whose union is the vertex set $V$. The abstraction of this identity leads to an object that we call a polynomial set map of binomial type. Such a structure, which may be thought of as a generating function satisfying an exponential-like property, has some of the same properties as the chromatic polynomial, yet its definition does not suggest any connection to graph theory.

The idea of studying the chromatic polynomial through a generating function such as
\[
\sum_{S\subseteq V}\chi_{G|_S}(x)\mathbf{S}
\]
appears to have first been considered by Tutte~\cite[Eq. (13)]{tutte}, who showed essentially that
\[
\sum_{S\subseteq V}\chi_{G|_S}(x)\mathbf{S}=
\left(\sum_{\substack{S\subseteq V\\\textrm{stable}}}\mathbf{S}\right)^x=
\left(\sum_{S\subseteq V}\chi_{G|_S}(1)\mathbf{S}\right)^x,
\]
where vertex sets are multiplied by the rule $\mathbf{S} \mathbf{T}=\mathbf{S}\mbox{\boldmath$\cup$}\mathbf{T}$ if $S\cap T=\emptyset$, and $\mathbf{S} \mathbf{T}=0$ otherwise. More recently, Lass~\cite{MR1861053} has applied this approach to enumeration of acyclic orientations.

We will show that, coupled with the notion of binomial type, the set map approach provides a useful tool for studying the chromatic polynomial. Our main result is that
\[
\chi_G(x)=\sum_{\sigma\vdash V}
a_{\ell(\sigma)}(x)
\prod_{T\in\sigma}A\chi_{G|_T}(x),
\]
where $a_n(x)$ is any polynomial sequence of binomial type (meaning it satisfies $a_n(x+y)=\sum_{k=0}^{n}\binom{n}{k}a_k(x)a_{n-k}(y)$), $A$ is the linear functional (on polynomials) defined by $Aa_n(x)=\delta_{n=1}$, and the sum is over all set partitions of $V$. This provides a description of the coefficients in the expansion of the chromatic polynomial in terms of any polynomial sequence of binomial type.

This result has several interesting consequences. For example, we will show that the coefficients in the expansion of the chromatic polynomial in terms of rising factorials $x^{(n)}:=x(x+1)(x+2)\cdots(x+n-1)$ can be given the following combinatorial interpretation, originally due to Brenti~\cite[Theorem 5.5]{MR1069745}. Let
\[
\chi_G(x)=\sum_{k\geq0}c_k x^{(k)}.
\]
Then $(-1)^{|V|-k}c_k$ is the number of pairs $(\sigma,\alpha)$, where $\sigma$ is a set partition of $V$ with $k$ blocks and $\alpha$ is an acyclic orientation of all edges of $G$ with both ends in the same block of $\sigma$.

Though the chromatic polynomial is our primary example of a polynomial set map of binomial type, we describe a less trivial example in Section~\ref{section:abel}. We call it the Abel polynomial set map because it generalizes the Abel polynomial sequence $x(x-an)^{n-1}$. The coefficients of Abel polynomials are well-known to enumerate planted forests, and we show that the coefficients of the Abel polynomial set map enumerate a similar, more general structure.

We begin in Section~\ref{section:setmaps} by studying the algebra of set maps. In Section~\ref{section:review} we briefly review some definitions and results from umbral calculus. This allows us to prove our main result (Theorem~\ref{thm:mix}) in Section~\ref{section:polysetmaps}. Section~\ref{section:chromatic} is devoted to applications to the chromatic polynomial. Finally, in Section~\ref{section:abel} we describe a set map generalization of Abel polynomials.

\section{Notation}

Let $\mathbb{K}$ be a field of characteristic zero. A polynomial sequence $a(x)$ is a sequence of polynomials $a_0(x),a_1(x),\ldots\in\mathbb{K}[x]$ such that $\mathrm{deg}(a_n(x))=n$. This restriction, that degrees agree with indices, is essential in umbral calculus because it ensures that polynomial sequences are bases for $\mathbb{K}[x]$. Because we will not need to regard polynomial set maps as bases, we do not make a similar restriction in defining polynomial set maps.

A set partition is a set of disjoint, nonempty sets. We use the symbols $\sigma$, $\pi$, $\gamma$, and $\rho$ exclusively for set partitions. The elements of a set partition are called blocks. The length $\ell(\sigma)$ of a set partition $\sigma$ is the number of blocks. If the union of the blocks is $S$, then we say that $\sigma$ is a set partition of $S$ and write $\sigma\vdash S$. The same notation is used for integer partitions, but the meaning should always be evident from context (we have set partitions of sets and integer partitions of integers). Note that the only set partition of $\emptyset$ is $\emptyset$.

We use the following notations:
\begin{itemize}
\item
$[n]=\{1,2,\ldots,n\}$.
\item
$(x)_n=x(x-1)(x-2)\cdots(x-n+1)=\frac{x!}{(x-n)!}$.
\item
$x^{(n)}=x(x+1)(x+2)\cdots(x+n-1)=\frac{(x+n-1)!}{(x-1)!}$.
\item
If $P$ is a proposition, then $\delta_P=\left\{\begin{array}{ll}
1 & \textrm{if $P$ is true}\\
0 & \textrm{if $P$ is false}\\
\end{array}\right..
$
\item
$T\uplus U=S$ means $T\cup U=S$ and $T\cap U=\emptyset$.
\end{itemize}

\section{Set Maps}\label{section:setmaps}
In this section we study the ring of set maps. This ring and numerous applications have been studied by Lass (see~\cite{MR1861053}). Our approach differs from Lass's in that we define composition directly by~\eqref{eq:comp}, whereas Lass defines it using generating functions. The latter approach has the advantage of trivializing most of the results of this section. However, our approach makes clear the relationship between composition and the combinatorics of set partitions.

Let $V$ be a set. If $p$ is a function whose domain is the set of all finite subsets of $V$ and whose range is some set $E$, we call $p$ an \emph{$E$-set map on $V$}, or simply an \emph{$E$-set map}. The image of $S\subseteq V$ will always be denoted by $p_S$. Throughout this section we let $R$ be a fixed commutative, associative ring.

\begin{defn}Let $g$ and $h$ be $R$-set maps. Define their sum $g+h$ and product $g\cdot h$ to be the $R$-set maps
\[
(g+h)_S=
g_S+h_S
\]
and
\[
(g\cdot h)_S=
\sum_{T\uplus U=S}g_T h_U.
\]
\end{defn}
It is clear that this makes the set of $R$-set maps into a commutative, associative ring with unit $\delta_{S=\emptyset}$. This is a generalization of the ring of exponential formal power series in the following sense. We can regard a sequence $a$ as the set map $a_S=a_{|S|}$. If $h$ is a set map and $h_S$ depends only on $|S|$, we can likewise regard $h$ as a sequence. If $a$ and $b$ are sequences, we have
\[
(a\cdot b)_n=
\sum_{k=0}^n \binom{n}{k}a_k b_{n-k},
\]
from which it follows that
\[
\mathrm{gen}_{a\cdot b}(t)=\mathrm{gen}_a(t)\mathrm{gen}_b(t),
\]
where for a sequence $c$, $\mathrm{gen}_c(t)$ denotes the exponential formal power series
\[
\mathrm{gen}_c(t)=
\sum_{k\geq0}\frac{1}{k!}c_k t^k.
\]
Hence the product of set maps generalizes the product of exponential formal power series. We next define the corresponding generalization of composition of exponential formal power series.

\begin{defn}Let $h$ be an $R$-set map with $h_\emptyset=0$. If $a$ is an $R$-sequence (a sequence of elements of $R$), define the composition of $a$ with $h$ to be the $R$-set map
\begin{equation}\label{eq:comp}
(a\circ h)_S=
\sum_{\sigma\vdash S}a_{\ell(\sigma)}
\prod_{T\in\sigma}h_T.
\end{equation}
\end{defn}

We adopt the convention that an empty product is $1$, so $(a\circ h)_\emptyset=a_0$. The \emph{compositional formula} (see~\cite[Thm. 5.1.4]{MR1676282} for a proof),
\[
\sum_{k\geq0}\frac{1}{k!}a_k
\left(
\sum_{j\geq1}\frac{1}{j!}b_j t^j
\right)^k
=
\sum_{k\geq0}\frac{1}{k!}t^k
\sum_{\sigma\vdash [k]}a_{\ell(\sigma)}
\prod_{T\in\sigma}b_{|T|},
\]
shows that for sequences $a$ and $b$ with $b_0=0$,
\begin{equation}\label{eq:gencirc}
\mathrm{gen}_{a\circ b}(t)=\mathrm{gen}_a(\mathrm{gen}_b(t)).
\end{equation}

When $V=\{v_1,v_2,\ldots,v_n\}$ is finite, we can regard an $R$-set map on $V$ as an element of the ring
\[
R[V]:=R[v_1,v_2,\ldots,v_n]/\langle v_1^2,v_2^2,\ldots,v_n^2 \rangle
\]
through the isomorphism
\[
h\mapsto\sum_{\substack{S\subseteq V}}h_S \mathbf{S},
\]
where $\mathbf{S}\in R[V]$ denotes the product of the elements of $S$. If $f\in R[V]$ has no constant term, we have $f^n/n!\in R[V]$ and $f^n=0$ for $n>|V|$, so it makes sense to compose an exponential formal power series with $f$. One can verify that products of elements of $R[V]$ and compositions of exponential formal power series with elements of $R[V]$ that have no constant term coincide with products of set maps and composition of sequences with set maps $h$ such that $h_\emptyset=0$. The ring $R[V]$ plays a central role in Lass's approach to set maps.

We will need the following lemma.
\begin{lem}\label{lem:unl}Let $a$ and $b$ be $R$-sequences, and let $h$ be an $R$-set map. Then
\[
(a\cdot b)\circ h=(a\circ h)\cdot(b\circ h).
\]
\end{lem}
\begin{pf}We have
\[
((a\cdot b)\circ h)_S=
\sum_{\sigma\vdash S}
(a\cdot b)_{\ell(\sigma)}
\prod_{W\in\sigma}h_W=
\sum_{\sigma\vdash S}
\sum_{\pi\uplus\gamma=\sigma}a_{\ell(\pi)} b_{\ell(\gamma)}
\prod_{W\in\sigma}h_W.
\]
Interchanging the order of summation yields
\[
\sum_{T\uplus U=S}
\left(
\sum_{\pi\vdash T}
a_{\ell(\pi)}
\prod_{W\in\pi}h_W
\right)
\left(
\sum_{\gamma\vdash U}
b_{\ell(\gamma)}
\prod_{W\in\gamma}h_W
\right)=
((a\circ h)\cdot(b\circ h))_S.
\]
\qed\end{pf}

We next consider multiplicative and compositional inverses.
\begin{prop}\label{prop:inverse}Let $h$ be an $R$-set map such that $h_\emptyset=1$. Then
\[
h^{-1}_S=\sum_{\sigma\vdash S}(-1)^{\ell(\sigma)}\ell(\sigma)!\prod_{W\in\sigma}h_W
\]
is the multiplicative inverse of $h$.
\end{prop}
\begin{pf}We have
\[
(h\cdot h^{-1})_S=\sum_{T\uplus U=S}h_T
\sum_{\sigma\vdash U}(-1)^{\ell(\sigma)}\ell(\sigma)!\prod_{W\in\sigma}h_W.
\]
When $S=\emptyset$, this is 1. Otherwise, split the sum into the cases where $T=\emptyset$ and $T\neq\emptyset$. The first addend is just $h^{-1}_S$. The second is
\[
\sum_{\substack{T\uplus U=S\\T\neq\emptyset}}h_T
\sum_{\sigma\vdash U}(-1)^{\ell(\sigma)}\ell(\sigma)!\prod_{W\in\sigma}h_W=
\sum_{\sigma\vdash S}\ell(\sigma)(-1)^{\ell(\sigma)-1}(\ell(\sigma)-1)!\prod_{W\in\sigma}h_W=
-h^{-1}_S
\]
because splitting a set into two blocks, the first nonempty, and choosing a set partition of the second is equivalent to choosing a set partition of the whole set and a distinguished block. Hence $(h\cdot h^{-1})_S=\delta_{S=\emptyset}$.
\qed\end{pf}

In $R[V]$, Proposition~\ref{prop:inverse} is equivalent to the obvious formula
\[
(1+F)^{-1}=\sum_{k\geq0}(-1)^k F^k,
\]
where
\[
F=\sum_{\substack{S\subseteq V\\S\neq\emptyset}}h_S \mathbf{S}.
\]

\begin{lem}\label{lem:useful1}Let $g$ be a $\mathbb{K}$-set map, and let $a$ be a $\mathbb{K}$-sequence with $a_0=g_\emptyset$ and $a_1\neq0$. Then there is a unique $\mathbb{K}$-set map $h$ with $h_\emptyset=0$ such that $a\circ h=g$.
\end{lem}
\begin{pf}
For $S\neq\emptyset$, write
\[
g_S=a_1h_S+\sum_{\substack{\sigma\vdash S\\\ell(\sigma)>1}}a_{\ell(\sigma)}\prod_{T\in\sigma}h_T.
\]
Solving for $h_S$, we obtain an expression for $h_S$ in terms of $h_T$ for $T\subset S$, so the result follows by induction.
\qed\end{pf}

\begin{lem}\label{lem:useful}Assume $V$ is countably infinite. Let $h$ and $g$ be $\mathbb{K}$-set maps with $h_\emptyset=0$. Let $a$ be a $\mathbb{K}$-sequence such that $a\circ h=g$. If $h$ is nonzero on one-element sets, then $a$ is uniquely determined by this equation.
\end{lem}
\begin{pf}
Let $n>0$, and let $S\subseteq V$ be such that $|S|=n$. Write
\[
g_S=a_n\prod_{v\in S}h_{\{v\}}+\sum_{\substack{\sigma\vdash S\\\ell(\sigma)<n}}a_{\ell(\sigma)}\prod_{T\in\sigma}h_T.
\]
Solving for $a_n$, we see that $a_n$ is determined by $a_1, \ldots, a_{n-1}$, so the result follows by induction.
\qed\end{pf}

\section{Umbral Calculus}\label{section:review}
In this section we review the concepts from umbral calculus that we will require. Proofs may be found in~\cite{MR741185} or~\cite{MR0485417}.

If $L,M\in\mathbb{K}[x]^*$ are linear functionals, their product $LM$ is defined to be the linear functional whose action on the monomials $x^n$ is given by
\[
LMx^n=\sum_{k=0}^n\binom{n}{k}Lx^k Mx^{n-k}.
\]
This product is clearly associative and commutative. With this product, the space $\mathbb{K}[x]^*$ is called the \emph{umbral algebra}.

A \emph{delta functional} $A\in\mathbb{K}[x]^*$ is a functional satisfying $Ax\neq0$ and $A1=0$. 
\begin{lem}\label{lem:unique}If $A$ is a delta functional and $f(x)$ is a polynomial, then $f(x)$ is uniquely determined by the sequence $A^k f(x)$ (where $A^0f(x)=f(0)$).
\end{lem}

A polynomial sequence $a(x)$ is said to be of \emph{binomial type} if
\begin{equation}\label{eq:stepumbral2}
a_n(x+y)=\sum_{k=0}^{n}\binom{n}{k}a_k(x)a_{n-k}(y)
\end{equation}
for all $n\geq0$. Notice that in the notation of the previous section,~\eqref{eq:stepumbral2} may be written as $a(x+y)=a(x)\cdot a(y)$.

\begin{thm}\label{thm:btdf}The following are equivalent.
\begin{enumerate}
\item
$a(x)$ is a polynomial sequence of binomial type.
\item
There is a delta functional $A$ satisfying $A^k a_n(x)=k!\delta_{n=k}$ for all $n,k\geq0$.
\end{enumerate}
\end{thm}

If $a(x)$ is a polynomial sequence of binomial type and $A$ is the delta functional of Theorem~\ref{thm:btdf}, then we say that $a(x)$ is \emph{associated} to $A$.

\section{Polynomial Set Maps}\label{section:polysetmaps}

If the range of a set map $p$ is the polynomial ring $\mathbb{K}[x]$, then we call $p=p(x)$ a \emph{polynomial set map}. We use the notation $(p(x))_S=p_S(x)$. The following definition appears in~\cite[Eq. 5.8]{MR544721} (with a different name), where it is attributed to J. P. S. Kung and T. Zaslavsky.
\begin{defn}Let $p(x)$ be a polynomial set map. If $p(x+y)=p(x)\cdot p(y)$, or equivalently,
\begin{equation}\label{eq:binset}
p_S(x+y)=\sum_{T\uplus U=S}p_T(x)p_U(y)
\end{equation}
for every finite subset $S\subseteq V$, then we say that $p(x)$ is of binomial type.
\end{defn}
If $p_S(x)$ depends only on $|S|$, \eqref{eq:binset} becomes
\[
p_n(x+y)=\sum_{k=0}^{n}\binom{n}{k}p_k(x)p_{n-k}(y),
\]
so if $\mathrm{deg}(p_n(x))=n$, then $p(x)$ is polynomial sequence of binomial type.

\begin{lem}\label{lem:mech}Let $p(x)$ be a polynomial set map of binomial type.
\begin{enumerate}
\item Either $p_\emptyset(x)=0$ or $p_\emptyset(x)=1$.
\item For every finite non-empty subset $S\subseteq V$, $p_S(0)=0$.
\end{enumerate}
\end{lem}
\begin{pf}We have $p_\emptyset(x+y)=p_\emptyset(x)p_\emptyset(y)$, so for every positive integer $t$ we have $p_\emptyset(t)=p_\emptyset(1)^t$. Hence the only way $p_\emptyset(x)$ can be a polynomial is if $p_\emptyset(x)=0$ or $1$. This proves the first part.

Now suppose $p_S(0)=0$ for $0<|S|<k$. If $|S|=k$, then
\[
p_S(0)=\sum_{T\uplus U=S}p_T(0)p_U(0)=p_S(0)p_\emptyset(0)+p_\emptyset(0)p_S(0)=2p_\emptyset(0)p_S(0),
\]
so $p_S(0)=0$ regardless of whether $p_\emptyset(x)=0$ or $1$. For the case $|S|=1$ we have
\[
p_S(0)=p_S(0)p_\emptyset(0)+p_\emptyset(0)p_S(0)=2p_\emptyset(0)p_S(0),
\]
so $p_S(0)=0$. Hence the second part follows by induction.
\qed\end{pf}

We say that $p(x)$ is \emph{nontrivial} if $p_\emptyset(x)=1$ and \emph{trivial} otherwise. If $p(x)$ is trivial, then
\[
p_S(x)=p_S(x+0)=\sum_{T\uplus U=S}p_T(x)p_U(0)=0,
\]
so there is only one trivial polynomial set map of binomial type.

\begin{lem}\label{lem:action}
If $p(x)$ is a polynomial set map of binomial type and $L,M\in\mathbb{K}[x]^*$ are linear functionals, then
$
LMp(x)=Lp(x)\cdot Mp(x).
$
\end{lem}
\begin{pf}The following proof is essentially the same as that of~\cite[Prop. 3.3]{MR0485417}. We denote by $L_x$ the linear operator on $\mathbb{K}[x,y]$ defined by $Lx^ny^m=(Lx^n)y^m$ (and similarly for $M_y$). The definition of the product $LM$ is then
\[
LMx^n=L_xM_y(x+y)^n.
\]
By linearity we may substitute any polynomial $f(x)$, thus obtaining
\[
LMf(x)=L_xM_yf(x+y).
\]
If $f(x)=p_S(x)$, we have
\[
LMp_S(x)=L_xM_yp_S(x+y)=
L_xM_y(p(x)\cdot p(y))_S=
(L_xp(x)\cdot M_yp(y))_S,
\]
from which the result follows.
\qed\end{pf}

We can now prove our main result.
\begin{thm}\label{thm:mix}Let $a(x)$ be a  polynomial sequence of binomial type associated to a delta functional $A$. If $p(x)$ is a nontrivial polynomial set map of binomial type, then
\[
p(x)=a(x)\circ Ap(x),
\]
or equivalently,
\[
p_S(x)=\sum_{\sigma\vdash S}
a_{\ell(\sigma)}(x)
\prod_{T\in\sigma}Ap_T(x)
\]
for every finite subset $S\subseteq V$.
\end{thm}
\begin{pf}
By Lemma~\ref{lem:unique} it will suffice to show that both sides are equal under the action of $A^k$. By Theorem~\ref{thm:btdf} the right hand side becomes
\begin{equation}\label{eq:aec}
\sum_{\substack{\sigma\vdash S\\\ell(\sigma)=k}}
k!
\prod_{T\in\sigma}Ap_T(x).
\end{equation}
Iterating Lemma~\ref{lem:action}, we see that the left hand side becomes
\[
\sum_{S_1\uplus S_2\uplus\cdots\uplus S_k=S}
~\prod_{i=1}^k A p_{S_k}(x)
\]
for $k>0$ and $\delta_{S=\emptyset}$ for $k=0$. Since $Ap_\emptyset(x)=A1=0$, this is actually a sum over ordered set partitions of $S$ with $k$ blocks. But there are $k!$ ways to order such a set partition, so this is the same as~\eqref{eq:aec}.
\qed\end{pf}

\begin{cor}\label{cor:bij}Let $a(x)$ be a  polynomial sequence of binomial type associated to a delta functional $A$. Then the map $\theta$ defined by $\theta(h)=a(x)\circ h$ is a bijection between $\mathbb{K}$-set maps $h$ such that $h_\emptyset=0$ and nontrivial polynomial set maps of binomial type. Furthermore, the inverse is given by $\theta^{-1}(p(x))=Ap(x)$.
\end{cor}
\begin{pf}We have
\[
a(x+y)\circ h=(a(x)\cdot a(y))\circ h=(a(x)\circ h) \cdot (a(y)\circ h)
\]
by Lemma~\ref{lem:unl}, which shows that $\theta(h)$ is indeed of binomial type. It is nontrivial because $(\theta(h))_\emptyset=a_0(x)=1$. Theorem~\ref{thm:mix} shows that $\theta\circ\theta^{-1}(p(x))=p(x)$, and it is clear from the definitions that $\theta^{-1}\circ\theta(h)=h$. Hence $\theta$ is a bijection.
\qed\end{pf}

\begin{cor}Assume $V$ is countably infinite. Let $p(x)$ be a polynomial set map of binomial type that is nonzero on one-element sets. Suppose that for some polynomial sequence $b(x)$ (not necessarily of binomial type) and function $F:\mathbb{K}[x]\rightarrow\mathbb{K}$ (not necessarily linear) we have
\[
p(x)=b(x)\circ Fp(x),
\]
where $(Fp(x))_S=F(p_S(x))$. Then $b(x)$ is a polynomial sequence of binomial type, and if we let $b(x)$ be associated to the delta functional $B$, then $Fp(x)=Bp(x)$.
\end{cor}
\begin{pf}
Since $p(x)$ is of binomial type, we have
\[
b(x+y)\circ Fp(x)=(b(x)\circ Fp(x))\cdot (b(y)\circ Fp(x))=(b(x)\cdot b(y))\circ Fp(x)
\]
by Lemma~\ref{lem:unl}. If $v\in V$, then $0\neq p_{\{v\}}(x)=b_1(x)F(p_{\{v\}}(x))$, so $F(p_{\{v\}}(x))\neq0$. It follows from Lemma~\ref{lem:useful} that $b(x+y)=b(x)\cdot b(y)$, so $b(x)$ is of binomial type. Let $b(x)$ be associated to the delta functional $B$. We have
\[
p(x)=b(x)\circ Bp(x)=b(x)\circ Fp(x).
\]
It follows from Lemma~\ref{lem:useful1} that $Bp(x)=Fp(x)$.
\qed\end{pf}

Theorem~\ref{thm:mix} makes it easy to expand a polynomial set map of binomial type in terms of any polynomial sequence of binomial type. Let us note some of the most useful examples.

For $a\in\mathbb{K}$, the Abel polynomial sequence $x(x-an)^{n-1}$ is associated to the delta functional $Lq(x)=q^\prime(a)$. Hence for any nontrivial polynomial set map of binomial type $p(x)$,
\begin{equation}\label{eq:te5}
p_S(x)=\sum_{\sigma\vdash S}x(x-a \ell(\sigma))^{\ell(\sigma)-1}
\prod_{T\in\sigma}p^\prime_T(a).
\end{equation}
For $a=0$, this becomes
\begin{equation}\label{eq:exp1}
p_S(x)=\sum_{\sigma\vdash S}x^{\ell(\sigma)}
\prod_{T\in\sigma}p^\prime_T(0).
\end{equation}
Evaluating at $x\in\mathbb{K}$, we obtain
\[
\sum_{S\subseteq V}p_S(x)\mathbf{S}=
\mathrm{exp}\left(x\sum_{S\subseteq V}p^\prime_S(0)\mathbf{S}\right)
\]
in $\mathbb{K}[V]$.

For $a\neq0$, the falling factorial polynomial sequence $(\frac{x}{a})_n$ is associated to the delta functional $Lq(x)=q(a)-q(0)$. By Lemma~\ref{lem:mech}, if $p(x)$ is a nontrivial polynomial set map of binomial type, then $p_S(0)=0$ for $S\neq\emptyset$. Hence
\begin{equation}\label{eq:te30}
p_S(x)=\sum_{\sigma\vdash S}(\frac{x}{a})_{\ell(\sigma)}
\prod_{T\in\sigma}p_T(a).
\end{equation}
Evaluating at $x\in\mathbb{K}$, we obtain
\[
\sum_{S\subseteq V}p_S(x)\mathbf{S}=
\left(\sum_{S\subseteq V}p_S(a)\mathbf{S}\right)^{x/a}
\]
in $\mathbb{K}[V]$ (recalling from Lemma~\ref{lem:mech} that $p_\emptyset(x)=1$), from which it follows that
\[
\sum_{S\subseteq V}p_S(xy)\mathbf{S}=
\left(\sum_{S\subseteq V}p_S(x)\mathbf{S}\right)^y.
\]

\section{The Chromatic Polynomial}\label{section:chromatic}
Our primary example of a polynomial set map of binomial type is the chromatic polynomial (see~\cite{MR95h:05105} for background). Traditionally, the chromatic polynomial is regarded as a function of a graph and a number (of colors). Here we will instead assume that $G$ is fixed and define the chromatic polynomial set map as a polynomial defined for each subset of the vertices of $G$. For a graph $G$ and a subset $S$ of the vertices, the \emph{restriction of $G$ to $S$} is the graph $G|_S$ with vertex set $S$ and edge set consisting of all edges of $G$ with both ends in $S$.

\begin{defn}Let $G$ be a finite graph with vertex set $V$. If $S\subseteq V$, define $\chi_S(x)=\chi_{G|_S}(x)$, where $\chi_{G|_S}(x)$ is the traditional chromatic polynomial of the restriction of $G$ to $S$. We call $\chi(x)$ the chromatic polynomial set map of $G$.
\end{defn}

The following fact (for the more general dichromatic polynomial) was observed by Tutte~\cite[Eq. (14)]{tutte}.

\begin{thm}\label{thm:chrombin}The chromatic  polynomial set map is of binomial type.
\end{thm}
\begin{pf}If we have a coloring of $G$ with $x+y$ colors, we can divide $V$ into two disjoint subsets depending on whether a vertex is colored with one of the first $x$ colors or with the remaining $y$ colors. Furthermore, if we restrict the coloring to one of the subsets, it will still be proper. Conversely, if we choose disjoint subsets and proper colorings of each, the first using only the first $x$ colors and the second using only the remaining $y$ colors, the induced coloring of $G$ will be proper. This proves the theorem for the case where $x$ and $y$ are positive integers. But the chromatic polynomial is a polynomial, so it must be true in general.
\qed\end{pf}

For the chromatic polynomial,~\eqref{eq:te5} and~\eqref{eq:te30} give the following expansions.
\begin{cor}For all $a\in\mathbb{K}$, we have
\begin{equation}\label{eq:cexp1}
\chi_S(x)=\sum_{\sigma\vdash S}x(x-a \ell(\sigma))^{\ell(\sigma)-1}
\prod_{T\in\sigma}\chi^\prime_T(a).
\end{equation}
\end{cor}
\begin{cor}For all nonzero $a\in\mathbb{K}$, we have
\begin{equation}\label{eq:cexp2}
\chi_S(x)=\sum_{\sigma\vdash S}(\frac{x}{a})_{\ell(\sigma)}
\prod_{T\in\sigma}\chi_T(a).
\end{equation}
\end{cor}

Setting $a=0$ in~\eqref{eq:cexp1}, we obtain the well-known expansion
\begin{equation}\label{eq:d0}
\chi_S(x)=
\sum_{\sigma\vdash S}x^{\ell(\sigma)}
\prod_{T\in\sigma}\chi^\prime_T(0).
\end{equation}
There are at least four known combinatorial interpretations of the number $(-1)^{|S|-1}\chi^\prime_S(0)$. It is the number of broken circuit free subtrees of $G|_S$~\cite[Sect. 7]{logicalexpansion}, the number of acyclic orientations of $G|_S$ with a unique sink at a fixed vertex~\cite[Thm. 7.3]{MR712251}, the number of cycle permutations of $S$ whose minimal expression as a product of transpositions contains only transpositions corresponding to the edges of $G|_S$~\cite{MR2048555}, and the number of increasing $G|_S$-connected trees~\cite{MR2134187}.

Since $\chi_S(1)$ is the number of proper colorings of $G|_S$ with just one color, $\chi_S(1)=1$ if $G|_S$ has no edges and $\chi_S(1)=0$ otherwise. Therefore, if we set $a=1$ in~\eqref{eq:cexp2}, we  obtain the classical expansion
\begin{equation}\label{eq:stab}
\chi_S(x)=\sum_{\substack{\sigma\vdash S\\\textrm{stable}}}(x)_{\ell(\sigma)},
\end{equation}
where a set partition $\sigma\vdash S$ is said to be stable if for each block $T\in\sigma$, $G|_T$ has no edges.

Setting $a=-1$ in~\eqref{eq:cexp2} yields
\[
\chi_S(x)=\sum_{\sigma\vdash S}(-1)^{\ell(\sigma)}x^{(\ell(\sigma))}
\prod_{T\in\sigma}\chi_T(-1).
\]
A result of Stanley~\cite[Cor. 1.3]{MR0317988} states that $(-1)^{|S|}\chi_S(-1)$ is the number of acyclic orientations of $G|_S$. Hence we have the following combinatorial interpretation of the coefficient of $x^{(k)}$ in $\chi_S(x)$. If $\sigma\vdash S$, we denote by $G|_\sigma$ the graph with vertex set $S$ and edge set consisting of all edges of $G$ with both ends in the same block of $\sigma$.
\begin{prop}[Brenti]\label{prop:exp91}Let
\[
\chi_S(x)=\sum_{k\geq0}c_k x^{(k)}.
\]
Then $(-1)^{|S|-k}c_k$ is the number of pairs $(\sigma,\alpha)$, where $\sigma\vdash S$, $\ell(\sigma)=k$, and $\alpha$ is an acyclic orientation of $G|_\sigma$.
\end{prop}
This result was first proved by Brenti~\cite[Theorem 5.5]{MR1069745}, with a simpler proof due to Lass~\cite[Proposition 6.1]{MR1861053}.

Setting $a=1$ in~\eqref{eq:cexp1} yields another interesting expansion.
\begin{prop}\label{prop:exp92}
\[
\chi_S(x)=
\sum_{\sigma\vdash S}x(x-\ell(\sigma))^{\ell(\sigma)-1}
\prod_{T\in\sigma}\chi^\prime_T(1).
\]
\end{prop}
The invariant $|\chi^\prime_S(1)|$ was introduced (for matroids) by Crapo~\cite{MR0215744}. By a theorem of Greene and Zaslavsky~\cite[Thm. 7.2]{MR712251}, if $G|_S$ has no isolated vertices and at least one edge, then $|\chi^\prime_S(1)|$ is the number of acyclic orientations of $G|_S$ with a unique sink and source at fixed adjacent vertices (see also~\cite{MR1778205} and~\cite{MR1823628}).

The polynomial sequence of binomial type
\[
b_n(x)=\sum_{\sigma\vdash[n]}(x)_{\ell(\sigma)}
\prod_{T\in\sigma}(-1)^{|T|-1}(|T|-1)!,
\]
which has the simple generating function
\[
\sum_{n\geq0}\frac{1}{n!}b_n(x)t^n=(1+\log(1+t))^x,
\]
is associated to the delta functional $B$ defined by $B(x)_n=1$ for $n>0$ and $B1=0$. The functional $B$ (except with $B1=1$) plays a prominent role in~\cite{MR0161805}. By~\eqref{eq:stab},
\[
B\chi_S(x)=\sum_{\substack{\sigma\vdash S\\\textrm{stable}}}1
\]
is the number of stable partitions of $S$. Hence we have the following expansion.
\begin{prop}\label{prop:exp93}Let $s_T$ be the number of stable partitions of $G|_T$. Then
\[
\chi_S(x)=
\sum_{\sigma\vdash S}b_{\ell(\sigma)}(x)
\prod_{T\in\sigma}s_T.
\]
\end{prop}
The coefficient of $b_n(x)$ in $\chi_S(x)$ is also the coefficient of $(x)_n$ in the $\sigma$-polynomial of $G|_S$ (as defined in~\cite{MR1069745}). Like the chromatic polynomial set map, the $\sigma$-polynomial set map is of binomial type.

To the knowledge of the author, Propositions~\ref{prop:exp92} and~\ref{prop:exp93} have not appeared in the literature.

\section{The Abel Polynomial Set Map}\label{section:abel}

In this section we will describe a nontrivial example of a polynomial set map of binomial type. If $\sigma$ is a set partition of a set $S$, define $\|\sigma\|=|S|$.

\begin{thm}\label{thm:abel}Let $\sigma\vdash V$. Define a polynomial set map $f(x)$ on $\sigma$ by
\[
f_\pi(x)=x(x+\|\pi\|)^{\ell(\pi)-1}
\]
for every finite subset $\pi\subseteq\sigma$. Then $f(x)$ is of binomial type.
\end{thm}

\psset{unit=2cm}
\psset{yunit=-2cm}
\begin{figure}
\begin{pspicture}(0,0)(7.0,6.0)

\psset{dotsize=0.06}
\psset{linewidth=0.04}
\psset{arrowscale=1.2}
\psset{origin={0,-2.85}}

\pscircle(1.0,1.0){0.42}
\pscircle(1.0,2.0){0.42}
\psframe*[linecolor=white](0.5,1.0)(1.5,2.0)
\psline(0.6,1.0)(0.6,2.0)
\psline(1.4,1.0)(1.4,2.0)

\pscircle(2.0,1.0){0.42}
\pscircle(2.0,1.5){0.42}
\psframe*[linecolor=white](1.5,1.0)(2.5,1.5)
\psline(1.6,1.0)(1.6,1.5)
\psline(2.4,1.0)(2.4,1.5)

\pscircle(2.0,2.5){0.42}
\pscircle(2.0,2.5){0.48}

\pscircle(3.0,1.0){0.42}
\pscircle(3.0,2.5){0.42}
\psframe*[linecolor=white](2.5,1.0)(3.5,2.5)
\psline(2.6,1.0)(2.6,2.5)
\psline(3.4,1.0)(3.4,2.5)

\pscircle(4.0,1.0){0.42}
\pscircle(4.0,2.0){0.42}
\psframe*[linecolor=white](3.5,1.0)(4.5,2.0)
\psline(3.6,1.0)(3.6,2.0)
\psline(4.4,1.0)(4.4,2.0)

\pscircle(5.0,1.0){0.42}
\pscircle(5.0,1.5){0.42}
\psframe*[linecolor=white](4.5,1.0)(5.5,1.5)
\psline(4.6,1.0)(4.6,1.5)
\psline(5.4,1.0)(5.4,1.5)

\pscircle(5.0,2.5){0.42}

\pscircle(6.0,1.0){0.42}
\pscircle(6.0,1.5){0.42}
\psframe*[linecolor=white](5.5,1.0)(6.5,1.5)
\psline(5.6,1.0)(5.6,1.5)
\psline(6.4,1.0)(6.4,1.5)

\psdots(1.0,1.0)
\psdots(1.0,1.5)
\psdots(1.0,2.0)

\psdots(2.0,1.0)
\psdots(2.0,1.5)

\psdots(3.0,1.0)
\psdots(3.0,1.5)
\psdots(3.0,2.0)
\psdots(3.0,2.5)

\psdots(4.0,1.0)
\psdots(4.0,1.5)
\psdots(4.0,2.0)

\psdots(5.0,1.0)
\psdots(5.0,1.5)

\psdots(6.0,1.0)
\psdots(6.0,1.5)

\psdots(2.0,2.5)
\psdots(5.0,2.5)

\psset{linewidth=0.03}

\psline{->}(1.4,2.0)(2.0,2.5)
\psline{->}(2.4,1.25)(3.0,1.0)
\psline{->}(3.4,1.75)(4.0,2.0)
\psline{->}(3.6,1.5)(3.0,1.5)
\psline{->}(5.0,1.9)(5.0,2.5)
\psline{->}(5.6,1.25)(5.0,1.5)
\psline{->}(5.4,2.5)(6.0,1.5)

\rput(1.0,0.4){$1$}
\rput(2.0,0.4){$2$}
\rput(3.0,0.4){$3$}
\rput(4.0,0.4){$4$}
\rput(5.0,0.4){$5$}
\rput(6.0,0.4){$6$}
\rput(2.0,3.15){$7$}
\rput(5.0,3.1){$8$}

\psset{origin={0,0}}

\psset{dotsize=0.06}
\psset{linewidth=0.04}
\psset{arrowscale=1.2}

\pscircle(1.0,1.0){0.42}
\pscircle(1.0,2.0){0.42}
\psframe*[linecolor=white](0.5,1.0)(1.5,2.0)
\psline(0.6,1.0)(0.6,2.0)
\psline(1.4,1.0)(1.4,2.0)

\pscircle(2.0,1.0){0.42}
\pscircle(2.0,1.5){0.42}
\psframe*[linecolor=white](1.5,1.0)(2.5,1.5)
\psline(1.6,1.0)(1.6,1.5)
\psline(2.4,1.0)(2.4,1.5)

\pscircle(2.0,2.5){0.42}
\pscircle(2.0,2.5){0.48}

\pscircle(3.0,1.0){0.42}
\pscircle(3.0,2.5){0.42}
\psframe*[linecolor=white](2.5,1.0)(3.5,2.5)
\psline(2.6,1.0)(2.6,2.5)
\psline(3.4,1.0)(3.4,2.5)

\pscircle(4.0,1.0){0.42}
\pscircle(4.0,2.0){0.42}
\psframe*[linecolor=white](3.5,1.0)(4.5,2.0)
\psline(3.6,1.0)(3.6,2.0)
\psline(4.4,1.0)(4.4,2.0)

\pscircle(5.0,1.0){0.42}
\pscircle(5.0,1.5){0.42}
\psframe*[linecolor=white](4.5,1.0)(5.5,1.5)
\psline(4.6,1.0)(4.6,1.5)
\psline(5.4,1.0)(5.4,1.5)

\pscircle(5.0,2.5){0.42}

\pscircle(6.0,1.0){0.42}
\pscircle(6.0,1.5){0.42}
\psframe*[linecolor=white](5.5,1.0)(6.5,1.5)
\psline(5.6,1.0)(5.6,1.5)
\psline(6.4,1.0)(6.4,1.5)

\psdots(1.0,1.0)
\psdots(1.0,1.5)
\psdots(1.0,2.0)

\psdots(2.0,1.0)
\psdots(2.0,1.5)

\psdots(3.0,1.0)
\psdots(3.0,1.5)
\psdots(3.0,2.0)
\psdots(3.0,2.5)

\psdots(4.0,1.0)
\psdots(4.0,1.5)
\psdots(4.0,2.0)

\psdots(5.0,1.0)
\psdots(5.0,1.5)

\psdots(6.0,1.0)
\psdots(6.0,1.5)

\psdots(2.0,2.5)
\psdots(5.0,2.5)

\psset{linewidth=0.03}

\psline{->}(1.4,2.0)(2.0,2.5)
\psline{->}(2.4,1.25)(3.0,1.0)
\psline{->}(3.4,2.5)(5.0,2.5)
\psline{->}(3.6,1.5)(3.0,1.5)
\psline{->}(5.0,2.1)(5.0,1.5)
\psline{->}(5.4,1.25)(6.0,1.5)
\psline{->}(2.4,2.5)(4.0,2.0)

\end{pspicture}
\caption{An example of the bijection used in the proof of Theorem \ref{thm:abel} for the case $k=1$ and $n=8$. On top is a tail tree with one block circled. Beneath is the corresponding collection of $7$ tails with no two originating at the same block and with no tail originating at the circled block.}\label{tailbijfig}
\end{figure}
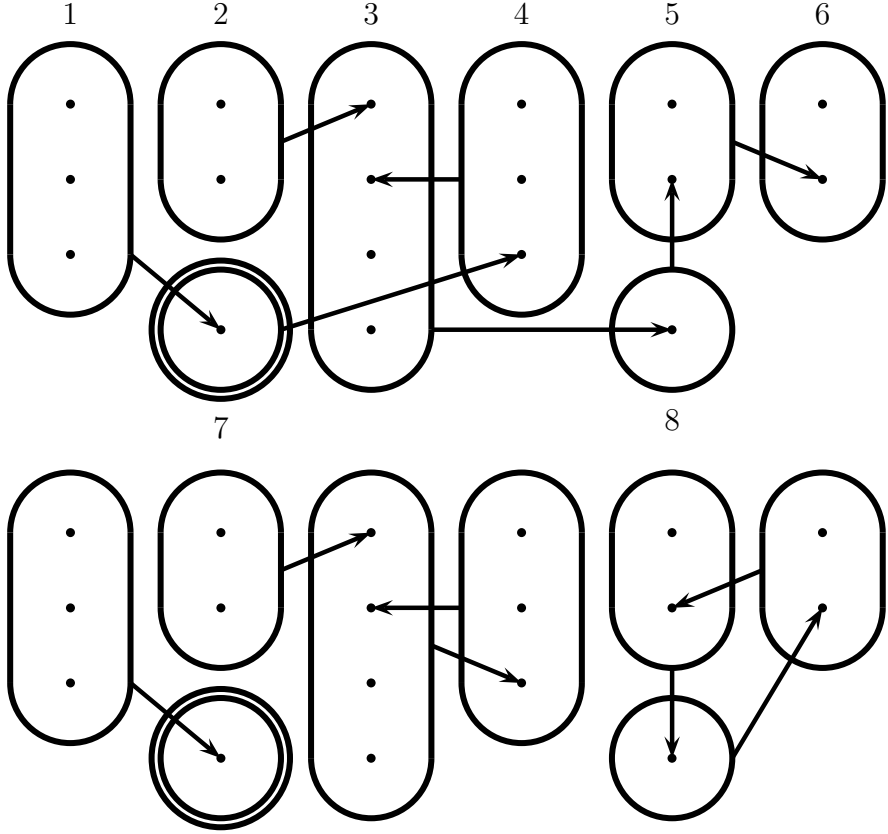

\begin{pf}We claim that
\begin{equation}\label{eq:sm3}
f_\pi(x)=\sum_{\gamma\vdash\pi}x^{\ell(\gamma)}
\prod_{\rho\in\gamma}\|\rho\|^{\ell(\rho)-1},
\end{equation}
which is of binomial type by Corollary~\ref{cor:bij}. Comparing coefficients of $x^k$ and letting $n=\ell(\pi)$, we see that this is equivalent to
\begin{equation}\label{eq:pfc}
\binom{n-1}{k-1}\|\pi\|^{n-k}=
\sum_{\substack{\gamma\vdash\pi\\\ell(\gamma)=k}}
\prod_{\rho\in\gamma}\|\rho\|^{\ell(\rho)-1}.
\end{equation}
If $\pi\vdash S$, define a \emph{tail on $\pi$}, or simply a \emph{tail}, to be an ordered pair $e=(T,v)\in \pi\times S$. We say that $e$ \emph{originates at} $T$ and \emph{points to} $v$, or that $v$ is the \emph{target} of $e$. See Figure \ref{tailbijfig} for a visualization of a collection of tails on a set partition. If $R$ is a set of tails, let $R/\pi$ be the digraph with vertex set $\pi$ obtained by changing each tail $(T,v)$ into a directed edge $(T,U)$, where $v\in U\in\pi$. Define a \emph{tail forest} to be a set $R$ of tails such that $R/\pi$ is acyclic and such that no two tails originate at the same block; that is, such that $R/\pi$ is a directed forest. We would like to show that the left hand side of~\eqref{eq:pfc} counts the number of tail forests with $k$ components, where the components of a tail forest $R$ are defined to correspond to the components of $R/\pi$. This would complete the proof, for it would show that $\|\rho\|^{\ell(\rho)-1}$ is the number of tail forests on $\rho$ with just one component (tail \emph{trees}), from which it would follow that the right hand side of~\eqref{eq:pfc} also counts the number of tail forests with $k$ components.

This can be proved by exhibiting a bijection between
\begin{itemize}
\item
pairs $(R,T)$, where $R$ is a tail forest with $k$ components and $T\in\pi$, and
\item
pairs $(F,T)$, where $F$ is a set of $n-k$ tails, no two of which originate at the same block, and $T\in\pi$ is a block at which no tail of $F$ originates.
\end{itemize}
We will describe such a map but omit the proof of bijectivity; it is similar to a standard enumeration of planted forests (see~\cite[Prop. 5.3.2, Second proof]{MR1676282}).

Let $\pi$ be a set partition with a fixed ordering of its blocks, and let $(R,T)$ be a pair as above. There is a unique root $r$ in the component of $R/\pi$ containing $T$, and a unique directed path from $T$ to $r$. Let $P$ be this path with the first block $T$ removed. For the tail tree shown in Figure \ref{tailbijfig}, $T=(4,3,8,5,6)$. This path can be regarded as the second line of a permutation written in two-line notation, with the first line being the elements of the path written in increasing order. This permutation can now be written in cycle notation; in our example, $(3,4)(5,8,6)$. Since we have a distinguished vertex in each block of $T$ given by the target of the tail originating at the previous block in the path, we can now replace the tails of the path with tails following the cycles of this permutation. Keeping the other tails of $R$ the same (except the tail originating at $T$, which must be removed), let $F$ be the resulting collection.

The construction of the inverse is now straightforward: form the cycles in $F$ into the corresponding path, and add $T$ to the beginning.
\qed\end{pf}

\begin{cor}\label{cor:abelsetmap}Suppose $\alpha$ is a $\mathbb{K}$-set map satisfying $\alpha_S=\sum_{v\in S}\alpha_{\{v\}}$. Then the polynomial set map
\[
p_S(x)=x(x+\alpha_S)^{|S|-1}
\]
is of binomial type.
\end{cor}
\begin{pf}
In the case where $\alpha_S$ is always a positive integer, this follows from Theorem~\ref{thm:abel} because there is a set partition whose block sizes correspond the numbers $\alpha_{\{v\}}$ for $v\in S$. Since we can regard $p_S(x)$ as a polynomial in the variables $\alpha_s$ for $s\in S$, it must also be true in general.
\qed\end{pf}

We call $p(x)$ the Abel polynomial set map because if $\alpha_S=-a|S|$, then $p(x)$ is the Abel polynomial sequence.

\section{Notes}\label{section:conclusion}

There are some important polynomial set maps of binomial type that we have not discussed. Graph theoretical examples include Tutte's dichromatic polynomial~\cite{tutte}, Zaslavsky's balanced chromatic polynomial of a signed graph~\cite{MR677061}, and the path and cycle polynomials of a digraph, whose product (in the set map sense) gives Chung and Graham's cover polynomial~\cite{MR1358990}. This latter polynomial set map has been investigated in detail by Lass~\cite{MR1928099}.

A useful generalization of polynomial sequences of binomial type is Sheffer sequences. Similarly, we can define a polynomial set map $f(x)$ to be a \emph{Sheffer set map} if there is a polynomial set map of binomial type $p(x)$ such that $f(x+y)=f(x)\cdot p(y)$. We say that $f(x)$ is \emph{nontrivial} if $f_\emptyset(x)\neq0$. An example of a Sheffer set map is $f_S(x)=\frac{1}{x}\chi_S(x)\delta_{v_0\in S}$, where $v_0\in V$ is a fixed vertex. The proof of Theorem~\ref{thm:mix} can be generalized to show that if $s(x)$ is a polynomial sequence that is Sheffer for $(N,A)$ (in the notation of~\cite[p. 17]{MR741185}), and if $f(x)$, where $f(x+y)=f(x)\cdot p(y)$, is a nontrivial Sheffer set map, then
\[
f(x)=Nf(x)\cdot (s(x)\circ Ap(x)),
\]
or equivalently,
\[
f_S(x)=\sum_{T\uplus U=S}Nf_T(x)
\sum_{\sigma\vdash U}
s_{\ell(\sigma)}(x)
\prod_{W\in\sigma}Ap_W(x)
\]
for every finite subset $S\subseteq V$.

Theorem~\ref{thm:mix} also admits a generalization where the concept of a polynomial set map of binomial type is replaced by the concept of a coalgebra map $C\rightarrow\mathbb{K}[x]$, where $C$ is a coalgebra and $\mathbb{K}[x]$ is given the coalgebra structure
\[
\Delta x^n=\sum_{k=0}^n\binom{n}{k}x^k\otimes x^{n-k}.
\]
This generalization can, for example, be used to obtain meaningful expansions of the order polynomial of a poset (see \cite[p. 320]{MR2037633}, \cite{MR0309813}, and \cite[p. 218]{MR1442260}).

Further generality can be obtained by replacing $\mathbb{K}[x]$ by a different coalgebra. For example, M\'endez's umbral calculus of symmetric functions \cite{MR1424311} can be used to generalize the results of Section \ref{section:chromatic} to Stanley's chromatic symmetric function \cite{MR96b:05174}.

\begin{ack}
I am grateful to Bruno Nachtergaele and Monica Vazirani for many helpful comments, and to Bodo Lass for pointing out to me the correct attribution of Proposition \ref{prop:exp91}.
\end{ack}

\bibliography{refs}

\end{document}